\documentclass{article}
\date{\today}
\usepackage{amsmath,amsfonts,color,amsthm}
\numberwithin{equation}{section} 

\newtheorem{lemma}[equation]{Lemma}
\newtheorem{proposition}[equation]{Proposition}

\newtheorem{remark}[equation]{Remark}
\newtheorem{definition}[equation]{Definition}
\renewcommand\i{{\boldsymbol{i}}}

\begin{document}
\title
{Some addition formulae for Abelian functions
for elliptic and hyperelliptic curves of
cyclotomic type}

\author{J. C. Eilbeck${}^{1}$, S. Matsutani${}^{2}$, Y. \^Onishi${}^{3}$\\
${}^{1}$Department of Mathematics and \\ the 
Maxwell Institute for Mathematical Sciences,\\
  Heriot-Watt University, Edinburgh EH14 4AS, UK\\ 
${}^{2}$ 8-21-1 Higashi-Linkan,\\ Sagamihara, 228-0811, Japan.\\
${}^{3}$ 4-3-11, Takeda, Kofu, 400-8511\\
Faculty of Education and Human Sciences\\
University of Yamanashi, Japan.
}

\label{firstpage}
\maketitle
\begin{abstract}{}
  We discuss a family of multi-term addition formulae for Weierstrass
  functions on specialized curves of low genus with many
  automorphisms, concentrating mostly on the case of genus one and
  two.  In the genus one case we give addition formulae for the
  equianharmonic and lemniscate cases, and in genus two we find some
  new addition formulae for a number of curves.
\end{abstract}

\section{Introduction}\label{sec1}

The aim of this paper is to introduce some new addition formulae for
the Weierstrass $\sigma$ and $\wp$ functions in genus one, and some
generalisations to some higher genus cases.  These formula are found
in the special case when some of the coefficients (moduli) of the
associated algebraic curves are chosen to be zero, and as a result the
curves have additional automorphisms (extra symmetries).
 
Although elliptic functions, including the Weierstrass elliptic
functions, have been extensively used (or perhaps over-used) to
enumerate travelling wave solutions of nonlinear wave equations,
relatively little has been written about the correspondingly higher
genus generalisations.  This is partly because no general handbooks
exist which play the same role as the familiar treatises on elliptic
functions.  This paper is part of a project to provide the material
for such a compendium.

Those coming to this paper because of possible applications to number
theory may prefer to see it as extending the classical theory of
complex multiplication for elliptic functions to higher genus
functions.  In this paper the complex multiplications are of
``cyclotomic type'', i.e. involving complex roots of unity.  These
generalise the more well-known addition formulae involving results for
$f(u+v)f(u-v)$ where we think of the $(\pm1)v$ as involving the two
real roots of unity.

We begin by summarising well-documented existing results for the genus
one case.  In this case we start with an elliptic curve reduced to the
standard Weierstrass form
\begin{equation}\label{curve_g1}
  y^2=4x^3-g_2x-g_3.
\end{equation}
The function $\wp(u)$  is the inverse function $u\mapsto x$ determined by
\begin{equation}
  u=\int_{\infty}^{(x,y)}\frac{dx}{2y},\label{pdef}
\end{equation}
and $\sigma(u)$ is an entire function satisfying
\begin{equation}\label{elliptic_sigma}
  \wp(u) = - \frac{d^2}{du^2}\log \sigma (u).
\end{equation}
The function $\wp(u)$ satisfies the well-known formula
\begin{equation}\label{WP}
  (\wp')^2=4\wp^3-g_2\wp-g_3,
\end{equation}
which is isomorphic to the genus one curve (but note this result does
not hold for the higher genus cases).  

The following two-variable addition formula plays an important role in
the theory of the Weierstrass $\sigma$ and $\wp$ functions, and its
generalisation are central to this paper:
\begin{equation}
  -\frac{\sigma(u+v)\sigma(u-v)}{\sigma(u)^2\sigma(v)^2}=\wp(u)-\wp(v).
   \label{A1}
\end{equation}
Taking the second logarithmic derivative of (\ref{A1}) gives the well
known addition formula involving just $\wp$ and $\wp'$, which is also
an addition formula on the curve (\ref{curve_g1}).

A three-variable addition formula is also known from the work of
Frobenius and Stickelburger \cite{fs77} (see also Whittaker \& Watson,
\cite{ww27})
\begin{equation}
\frac{\sigma(u-w)\sigma(v-w)\sigma(u-v)\sigma(u+v+w)}
{\sigma(u)^3\sigma(v)^3\sigma(w)^3} = -\frac12\left|
  \begin{array}{rrr}
   1 & \wp(u) & \wp'(u)\\
   1 & \wp(v) & \wp'(v)\\
   1 & \wp(w) & \wp'(w)
  \end{array}
\right|. \label{FS}
\end{equation}

In the genus two case, starting with the hyperelliptic curve 
\begin{equation}
  y^2=x^5+\mu_2 x^4+\mu_4 x^3+\mu_6 x^2+\mu_8 x +\mu_{10}, \label{C2}
\end{equation}
one can define generalized $\sigma$ and $\wp$ functions (see the
classical book by Baker \cite{ba07}, or Buchstaber {\em et al.}
\cite{bel97} for a modern treatment).  The main difference is that
$\sigma$ is now a function of $g=2$ variables, $u=\{u_1,u_2\}$, and
there are now three possible versions of the $\wp$ function, due to
the different possible logarithmic differentials of the $\sigma$
function:
\begin{equation}
  \wp_{ij}(u) = - \frac{\partial^2}{\partial u_i\partial u_j}\log \sigma (u), 
  \quad 1 \leq i\leq j \leq 2.
\end{equation}
(Note that in this notation the $\wp$ of the genus one theory would be
written as $\wp_{11}$).  The functions $\wp_{ij}$ and $\wp_{ijk}=
\frac{\partial}{\partial u_k} \wp_{ij}(u)$ satisfy equations analogous
to (\ref{WP}) The genus two $\sigma$ and $\wp$ functions satisfy an
analogue of the elliptic addition formula (\ref{A1})
\begin{equation}
  -\frac{\sigma(u+v)\sigma(u-v)}{\sigma(u)^2\sigma(v)^2}=
  \wp_{11}(u)-\wp_{11}(v)+\wp_{12}(u)\wp_{22}(v)-\wp_{22}(u)\wp_{12}(v)
   \label{A2}
\end{equation}
as described in Baker \cite{ba98,ba07}.  A generalisation of (\ref{FS}) in
the genus two case has also been derived (Eilbeck {\em et al.} \cite{eep03}),
but this is a rather complicated formula.

Our main aim in this paper is to point out that Abelian functions
associated with a curve with many automorphisms, namely, with extra
symmetries relative to general case, have novel addition formulae, which
are {\em not} valid in the general case.  In addition, we wish to
present the addition formulae in a way which makes the extra
symmetries explicit.

Although we restrict ourselves mostly to elliptic ($g=1$) and
hyperelliptic curves ($g=2$) in this paper, we comment briefly on
similar results which have been derived or are under study for more
general curves.

As an example, whilst the only nontrivial automorphism of the curve
(\ref{curve_g1}) with generic values of $g_j$s is $(x,y)\mapsto (x,-y)$,
the special curve
\[
 y^2=x^3-g_3 \ \ \mbox{with $g_3\neq 0$}
 \]
 has six automorphisms $(x,y)\mapsto(\zeta^j x, \pm y)$ with $j=0$,
 $1$, $2$ and $\zeta=\exp(2\pi i/3)$.  This has other addition
 formulae different from (\ref{A1}).  We mention here one from Theorem
 \ref{Te1}:
\[
   -\frac{\sigma(u \pm v)\sigma(u \pm \zeta v)\sigma(u \pm \zeta^2 v)}
{\sigma(u)^3\sigma(v)^3}=\pm \frac12\left(\wp'(u)\pm \wp'(v)\right).
\]
The formula is novel in the sense that, although it can be derived
from (\ref{A1}), it is only valid in case that $g_2=0$.

As another example, while the curve (\ref{C2}) with generic parameters 
has only two automorphisms $(x,y)\mapsto (x,\pm y)$, the special curve 
\begin{equation}
  y^2=x^5+\mu_{10},\quad\mbox{with $\mu_{10}\neq 0$},
\end{equation}
has ten automorphisms, $(x,y)\mapsto(\zeta^i x,\pm y),j=0,1,\dots,4$, with
$\zeta=\exp(2\pi i/5)$, and Abelian functions on its Jacobian variety
has addition formula different from (\ref{A2}), 
for example Proposition \ref{equipentamic} that expresses
\begin{equation}\label{LHS}
    \frac{ \sigma(u+v)\sigma(u+[\zeta]v)\sigma(u+[\zeta^2]v)
      \sigma(u+[\zeta^3]v)\sigma(u+[\zeta^4]v)}
    {\sigma(u)^5\sigma(v)^5},
\end{equation}
as a polynomial in $\wp_{ij}(u)$, $\wp_{ij}(v)$, $\wp_{ijk}(u)$, and
$\wp_{ijk}(v)$. Here $\zeta=\exp(2\pi i/5)$ and
$[\zeta^j]v=[\zeta^j](v_1,v_2) = (\zeta^j v_1,\zeta^{2j} v_2)$.

In genus one case, we obtain also three-term and four-term addition
formulae by using (\ref{FS}). While formulae of a similar type exist in
the higher genus cases, we do not mention them here because of their
complexity.

We have two ways to prove these formulae.  One is by simplifying the
expression in terms of $\wp$ and its derivative given by taking
product of modified formulae from (\ref{A1}) or (\ref{FS}), or their
higher genus generalizations.  The other is similar to the method used
in \cite{eemop07}, that is by balancing an expression such as
(\ref{LHS}) with a linear combination of suitable $\wp$-functions, in
which the derivation of the correct coefficients are aided by
algebraic computing software.

The paper is laid out as follows.  We first cover some basic theory
mainly needed for genus two and higher.  In \S \ref{s1}, we review the
types of curves we consider in the paper.  After introducing basic
notions in \S \ref{s2}, we define the function $\sigma(u)$ in \S
\ref{s3}, and the $\wp$-functions in \S \ref{s4}.  We consider the
genus one case in \S \ref{s5}, giving some detail to provide a
pedagogical background for the general methods.  The genus two case is
discussed in \S \ref{s6}.  In \S \ref{s7} we discuss briefly further
generalisations to higher genus cases, a topic which will be covered
in more detail elsewhere.

\section{Elliptic and hyperelliptic curves of cyclotomic type}
\label{s1}
In this section, we describe clearly the curves which we shall consider. 
Let $a\geq 2$ and $m$ be positive integers.  
We consider two type of curves according to \(am\) is odd or even. 
Namely, let
  \begin{equation}\label{def_eq}
  f(x,y)=
  \left\{
  \begin{aligned}
  &\begin{aligned}
  &y^2+\mu_{am}\,y-(x^{am}+\mu_{2a}x^{a(m-1)}+\mu_{4a}x^{a(m-2)}+\cdots\\
  &   \qquad\qquad\qquad\qquad +\mu_{2a(m-1)}x^a+\mu_{2am}), 
  \qquad\ \ \mbox{if $am$ is odd},\\
  \end{aligned}\\
  &\begin{aligned}
  &y^2-(x^{am+1}+\mu_{2a}x^{a(m-1)+1}+\mu_{4a}x^{a(m-2)+1}+\cdots\\
  & \qquad\qquad\qquad\qquad +\mu_{2a(m-1)}x^{a+1}+\mu_{2am}x),
  \ \ \ \mbox{if $am$ is even}.
  \end{aligned}
  \end{aligned}
  \right.
  \end{equation}
We consider the projective curve $C$ defined by the affine equation
  \[
  f(x,y)=0
  \]
with adding the unique point \(\infty\) at infinity. 
The genus of  $C$  is  
\[
  g=\lfloor am/2\rfloor
\] 
if it is non-singular. 
We refer the curve that is defined by the former equation as the
$(2,a[m])$-curve, and the later as the $(2,a[m]{+}1)$-curve.  Here the
first entry the number \lq\lq$2$" indicates these curves are either elliptic or
hyperelliptic curves, namely, the power of \(y\) in the defining equation.  
We are aiming to treat any algebraic curves
with a unique point at infinity, but we restrict ourselves here to
elliptic and hyperelliptic curves and to present our idea simply.

Both the $(2,a[m])$-curve and the $(2,a[m]{+}1)$-curve are acted on by 
the group $W_{2a}$ of $2a$-th roots of $1$ as automorphisms:
  \begin{equation}\label{action}
  \begin{aligned}
    {[-\zeta^2]} : (x,y)\mapsto& (\zeta^2 x,-y-\mu_{am}) & \mbox{ for a $(2,a[m])$-curve,}\\
    {[-\zeta]} : (x,y)\mapsto& (-\zeta x,{\i}y) & \mbox{ for a
      $(2,a[m]{+}1)$-curve,}
  \end{aligned}
  \end{equation}
where $\zeta$ is an $2a$-th root of $1$, and ${\i}^2=-1$. 

\noindent
{\bf Examples:} We give some examples here:

\begin{enumerate} 
\setlength{\itemsep}{1.0pt}
\item The general $(2,3[1])$-curve is defined by $y^2+\mu_3y=x^3+\mu_6$. 
  This is acted on by $W_6$.
\item The general $(2,2[1]{+}1)$-curve is defined by
  $y^2=x^3+\mu_4x$. This is acted on by $W_4$.
\item The general $(2,2[2]{+}1)$-curve is defined by
  $y^2=x^5+\mu_4x^3+\mu_8x$, which is the famous Burnside curve. 
This is acted on by $W_4$. 
\item The general $(2,5[1])$-curve is defined by
  $y^2+\mu_5y=x^5+\mu_{10}$. This is acted on by $W_{10}$.
\item The general $(2,4[1]{+}1)$-curve is defined by
  $y^2=x^5+\mu_8x$. This is acted on by $W_8$.
\item The general $(2,3[2]{+}1)$-curve is defined by
  $y^2=x^7+\mu_6x^4+\mu_{12}x$. This is acted on by $W_6$.
\item The general $(2,2[3]{+}1)$-curve is defined by
  $y^2=x^7+\mu_4x^5+\mu_8x^3+\mu_{12}x$. This is acted on by $W_4$.
\item The general $(2,7[1])$-curve is defined by
  $y^2+\mu_7y=x^7+\mu_{12}x$. This is acted on by $W_{14}$.
\end{enumerate}

In this paper we suppose that $\mu_j$ is $0$, if it does not appear in
the equation of $C$, $f(x,y)=0$.

\section{Differential forms,  etc.}\label{s2}
In the hyperelliptic case, the space of differential forms are spanned
by
\begin{equation}
   \omega_1=\frac{ dx}{2y},  \quad
   \omega_2=\frac{xdx}{2y},  \quad
   \cdots, \quad
   \omega_g=\frac{x^{g-1}dx}{2y},
   \label{eq1.02} 
\end{equation}
For variable $g$ points $(x_1,y_1)$, $(x_2,y_2)$, $\cdots$,
$(x_g,y_g)$ on $C$, we consider the integrals
\begin{eqnarray}
  u &=&(u_1,u_2,\cdots,u_g) \nonumber\\
  & =&\int_{\infty}^{(x_1,y_1)}\omega+\int_{\infty}^{(x_2,y_2)}\omega+
  \cdots  +\int_{\infty}^{(x_g,y_g)}\omega,
 \label{eq1.03}
\end{eqnarray}
where
\begin{equation}
  \omega=(\omega_1,\omega_2,\cdots,\omega_g).
  \label{eq1.04}
\end{equation}
Let
\begin{equation}
  \eta_j
  =\frac{1}{2y}\sum_{k=j}^{2g-j}(k+1-j) \mu_{4g-2k-2j} x^k dx
  \quad (j=1, \cdots, g), 
\end{equation}
which are differential forms of the second kind without poles 
except at  $\infty$.  

Let  $\Lambda$  be the lattice in  $\mathbb{C}^g$  
generated by the loop integrals of $\omega$\,: 
\begin{equation}\label{eq1.05}
  \Lambda=\bigg\{\oint \omega\bigg\}.
\end{equation}
Then the Jacobian variety of  $C$  is given by  $\mathbb{C}^g/\Lambda$. 
For $k=1$, $2$, $\cdots$, $g$, the map
\begin{eqnarray}\label{eq1.06}
 \iota : \mbox{Sym}^k(C)&\rightarrow& J, \nonumber\\
  (P_1,\cdots,P_k) &\mapsto& \left(\int_{\infty}^{P_1}\omega+\cdots+
  \int_{\infty}^{P_k}\omega\right) \mbox{ mod}\,\Lambda,
\end{eqnarray}
is an injection outside a certain small dimensional (relative to $k$)
subset.  If $k=g$, the map is surjective.  We denote the image
$\iota(\mathrm{Sym}^{k}(C))$ by $\Theta^{[k]}$.  Let
\newcommand{\rslt}{\mathrm{rslt}}
\begin{eqnarray}
R_1&=&\rslt_x\big(\rslt_y\big(f(x,y), f_x(x,y)\big),
                 \rslt_y\big(f(x,y), f_y(x,y)\big)\big), \nonumber\\
R_2&=&\rslt_y\big(\rslt_x\big(f(x,y), f_x(x,y)\big), 
                 \rslt_x\big(f(x,y), f_y(x,y)\big)\big), \\
R_3&=&\gcd(R_1,R_2),\nonumber
\end{eqnarray}
where $\rslt_z$ represents the resultant, namely, the determinant of
the Sylvester matrix with respect to the variable $z$.
Then $R_3$ is a perfect square 
in the ring
\[
  \mathbb{Z}[\{\mu_{j}\}]. 
\]
Hence we define 
\begin{equation}\label{discriminant}
  D={R_3}^{1/2}.
\end{equation}

\section{The sigma function}\label{s3}
\subsection{The definition of  $\sigma(u)$} 
We define here an entire function $\sigma(u)=\sigma(u_1,\cdots,u_g)$
on $\mathbb{C}^g$ associated with $C$, which we call the {\it  $\sigma$-function}.  
As usual, let
  \begin{equation}\label{sigma02}
  \alpha_i,\  \beta_j  \ \ \ (1\leq i,\ j \leq g)
  \end{equation} 
be closed paths on $C$ which generate $H_1(C,\mathbb{Z})$ such that
their intersection numbers are
$\alpha_i\cdot\alpha_j=\beta_i\cdot\beta_j=0$,
$\alpha_i\cdot\beta_j=\delta_{ij}$.

Define the period matrices by
\begin{equation}
  \left[\,\omega'  \ \omega''  \right]= 
  \left[\int_{\alpha_i}\omega_j \quad 
    \int_{\beta_i}\omega_j\right]_{i,j=1,\cdots,g},
  \,\,
  \left[\,\eta'  \ \eta''  \right]= 
  \left[\int_{\alpha_i}\eta_j \quad 
    \int_{\beta_i}\eta_j\right]_{i,j=1,\cdots,g}.
   \label{sigma03}
\end{equation} 
From (\ref{eq1.02}) we see the canonical divisor class of $C$ is given
by $4\infty$, and we are taking $\infty$ as the base point of the Abel
map (\ref{eq1.06}) for $k=g$. Hence the Riemann constant is an element
of $\big(\frac12\mathbb{Z}\big)^{2g}$ (see Mumford \cite{mu83}),
Coroll.3.11, p.166). Let
\begin{equation}
   \delta=\bigg[\begin{array}{cc}\delta'\\ \delta''\end{array}\bigg]
   \in \big(\tfrac12\mathbb{Z}\big)^{2g}
   \label{eq2.9} 
\end{equation} 
be the theta characteristic which gives the Riemann constant with
respect to the base point $\infty$ and the period matrix $[\omega'\
\omega'']$.  Note that we use $\delta',\delta''$ as well as $n$ in
(\ref{def-sigma}) as columns, to keep the notation a bit simpler.  
We define
  \begin{equation}\label{def-sigma}
  \begin{aligned}
    \sigma(u)&=\sigma(u_1,\cdots,u_g)\\
    &=c\,\exp(-\frac{1}{2}u\eta'{\omega'}^{-1}{}^t\negthinspace u)
    \Theta\negthinspace
    \left[\delta\right]({\omega'}^{-1}\ ^t\negthinspace u;\ {\omega'}^{-1}\omega'')\\
    &=c\,\mbox{exp}(-\frac{1}{2}u\eta'{\omega'}^{-1}\ ^t\negthinspace
    u) \sum_{n \in \mathbb{Z}^g}
    \exp\Big[2\pi \i\Big\{\frac12{}^t\negthinspace (n+\delta'){\omega'}^{-1}\omega''(n+\delta')\\
    &\qquad\qquad\qquad\qquad+{}^t\negthinspace
    (n+\delta')({\omega'}^{-1}\,^tu+\delta'')\Big\}\Big],
  \end{aligned}
  \end{equation}
where 
\begin{equation}
    c=\frac{1}{\sqrt[8]{D}}\bigg(\frac{\pi^g}{|\omega'|}\bigg)^{1/2}
   \label{sigma-const}
\end{equation}
with $D$ from (\ref{discriminant}).  Here the sign of a root of
(\ref{sigma-const}) is chosen so that the leading terms of $\sigma(u)$ is
just the Schur-Weierstrass polynomial that is originally defined in
Buchstaber {\em et al.} \cite{bel99}).
However, we use the modified version from \^Onishi \cite{on05b},
pg.~711.  The series (\ref{def-sigma}) converges because the imaginary
part of ${\omega'}^{-1}\omega''$ is positive-definite.

In what follows, for a given $u\in\mathbb{C}^g$, we denote by $u'$ and
$u''$ the unique elements in $\mathbb{R}^g$ such that
\begin{equation}
  u=u'\omega'+u''\omega''.
   \label{eq2.12}
\end{equation}
Then for $u$, $v\in\mathbb{C}^g$, and $\ell$
($=\ell'\omega'+\ell''\omega''$) $\in\Lambda$, we define
\begin{eqnarray}
  L(u,v)    &:=&{u}(\eta'{}^tv'+\eta''{}^tv''),\nonumber \\
  \chi(\ell)&:=&\exp[\pi \i \big(2({\ell'}\delta''-
  {\ell''}\delta') +{\ell'}{}^t\ell''\big)] \ (\in \{1,\,-1\}).
   \label{eq2.13}
\end{eqnarray}
In this situation, the most important properties of $\sigma(u;M)$
are as follows: 
\begin{lemma}\label{L2.14}
The function $\sigma(u)$ is an entire function which 
is independent of choice of the paths \(\alpha_j\), \(\beta_j\) of {\rm (\ref{sigma02})}. 
For all  $u\in\mathbb{C}^g$, $\ell\in\Lambda$ and we have 
  \begin{align}
  \sigma(u+\ell)&=\chi(\ell)\sigma(u)\exp L\Big(u+\frac12\ell,\ell\Big),\label{L2.14.1}\\
  u\mapsto\sigma(u)&\ \mbox{has zeroes of order $1$ along \(\Theta^{[g-1]}\)}, \label{L2.14.3}\\
  \sigma(u)&=0 \iff u\in\Theta^{[g-1]},\label{L2.14.4}
  \end{align}
where $\Theta^{[g-1]}$  is as defined following {\rm (\ref{eq1.06})}.
\end{lemma}
\begin{proof}
These are essentially classical results, and can be proved 
as in \cite{eemop07}, 
Lemma 4.1.  So we omit the proof. 
\end{proof}
\begin{lemma}\label{L2.15}
The coefficients in the expansion of the function \(\sigma(u)\) 
at the origin are polynomials of \(\mu_j\)s in {\rm (\ref{def_eq})} 
over the rationals \(\mathbb{Q}\). 
\begin{proof}
See Nakayashiki \cite{na08}. 
\end{proof}
\end{lemma}
\begin{lemma}\label{one-dim}
Let $\chi$ and $L$ be defined as above.  The space of entire
functions $\varphi(u)$ on $\mathbb{C}^g$ satisfying
  \[
  \varphi(u+\ell)=\chi(\ell)\varphi(u)\exp L\Big(u+\frac12\ell,\ell\Big)
  \]
is $1$-dimensional. 
\end{lemma}

\begin{proof}
  This is shown by the fact that the Pfaffian of the Riemann form
  attached to $L(\ ,\ )$ is $1$ (see Lang \cite{la82}, 
p.93, Th.3.1).
\end{proof}

\begin{definition}\label{weight}
  We introduce a \,{\rm weight}\, by letting the weight of $\mu_j$ to be
  $-j$, that of $u_j$ to be $2(g-j)+1$, that of $x$ to be $-2$, and that
  of $y$ to be $-2g-1$.
\end{definition}

We can easily show that any formula on $J$ is a sum of terms
homogeneous in this weight.  In many cases, the computation will be
easier if the weight is taken into consideration.  We can further
simplify by subdividing the calculations according to the separate
weights of the $\mu_i$ and the $u_j$ terms.

\subsection{Complex multiplication of $\sigma(u)$}

\begin{lemma}\label{CM}
  Let $C$ be a $(2,a[m])$- or $(2,a[m]{+}1)$-curve.  Let $\sigma(u)$
  be the sigma function associated with $C$ as above.  Let
\begin{eqnarray*}
  w&=& \left\{\begin{array}{ll}
  ((am)^2-1)/8, & \mbox{if $am$  is odd},\\
   ((am+1)^2-1)/8, & \mbox{if $am$ is even}.
   \end{array} \right.
\end{eqnarray*}
Let  $\zeta=\exp(2\pi\i/(2a))$. 
By the map \ref{eq1.06}, 
the action  (\ref{action}) on the curve \(C\) induce naturally
an action of the group  \(W_{2a}\) of \(2a\)-th roots of \(1\) on 
the space \(\mathbb{C}^g\) where \(u\) varies. 
We write this action explicitly for each cases below. 
Then we see that
  \[
  \sigma([-\zeta]u)=(-\zeta)^w\sigma(u).
  \]
\end{lemma}

\begin{proof}
Let  $\Lambda$  be the lattice in $\mathbb{C}^g$ as above. 
Then we have
  \[
  [-\zeta]\Lambda=\Lambda.
  \]
By (\ref{one-dim}), there is a constant $K$ such that
  \[
  \sigma([-\zeta]u)=K\sigma(u).
  \] 
Because $[-\zeta]^{2a}$ is the identity on  $\mathbb{C}^g$, we see
  \[
  K^{2a}=1.
  \]
By looking at the leading terms of $\sigma(u)$, we have
  \[
  K=(-\zeta)^w,
  \]
as desired. \end{proof}
\section{$\wp$-functions}\label{s4}
Using the sigma functions defined in the previous section, we let  
\begin{equation}\label{wp_fcts}
\wp_{jk}(u)=-\frac{\partial^2}{\partial u_j\partial u_k}\log\sigma(u), \quad
\wp_{jk\ell}(u)=\frac{\partial}{\partial u_\ell}\wp_{jk}(u), \quad \mbox{etc.}
\end{equation}
Then by (\ref{L2.14.1}), these functions are periodic with respect to
the $\Lambda$ of (\ref{eq1.05}).  If the genus of $C$ is $g=1$, then,
as usual, we write more classically $\wp_{11}(u)=\wp(u)$ and
$\wp_{111}(u)=\wp'(u)$.
\section{Genus One}
\label{s5}

\subsection{Generalities}

For completeness we start off with $C$ be the general elliptic curve
defined by
\[
  y^2+(\mu_1x+\mu_3)y=x^3+\mu_2x^2+\mu_4x+\mu_6.
\]
Then the $\wp(u)$ defined by (\ref{wp_fcts}) satisfies
  \begin{equation*}
  \wp'(u)= 2y+\mu_1x+\mu_3, \ \ \wp(u)=x
  \end{equation*}
when 
  \begin{equation*}
  u=\int_{\infty}^{(x,y)}\frac{dx}{2y+\mu_1x+\mu_3}, 
  \end{equation*}
and the  $\sigma(u)$, $\wp(u)$  satisfy (\ref{A1}) in the Introduction.

\subsection{Equianharmonic case}
\label{s5e}

We now specialize $C$ to the curve $y^2+\mu_3y=x^3+\mu_6$.  Then we have
\begin{equation}\label{2.02}
  (\wp')^2 = 4 \wp^3 + 4({\mu_3}^2+\mu_6). 
\end{equation}
As usual, by putting  $g_3=-4({\mu_3}^2+\mu_6)$, we rewrite (\ref{2.02}) as
\begin{equation}\label{2.04}
  (\wp')^2 = 4 \wp^3 - g_3. 
\end{equation}
This is usually called the {\em equianharmonic} case (see Abramowitz \&
Stegun \cite{as72}).  
Let $\zeta=\exp(2\pi\i/3)$. Then $\zeta^2 = -\zeta -1$,
and
\begin{equation}\label{CM1a}
\sigma(\zeta u) = \zeta \sigma(u), \quad \wp(\zeta u) = \zeta \wp(u), 
\quad \wp'(\zeta u) = \wp'(u),  
\end{equation}
by (\ref{CM}). 

The main results for the equianharmonic case are two novel addition
formulae, one for two variables and one for three variables:
\begin{proposition}\label{Te1}
\begin{align}\label{A1e2}
& -\frac{\sigma(u \pm v)\sigma(u \pm \zeta v)\sigma(u \pm \zeta^2 v)}
{\sigma(u)^3\sigma(v)^3}=\pm \frac12\left(\wp'(u)\pm \wp'(v)\right).\\
\label{A1e3}
& \frac{\sigma(u+v+w)\sigma(u+\zeta v+\zeta^2 w )
\sigma(u+\zeta^2 v+\zeta w )}{\sigma(u)^3\sigma(v)^3\sigma(w)^3}=\\
&\qquad  \frac14 (\wp'(u)\wp'(v)+\wp'(u)\wp'(w)+\wp'(v)\wp'(w))
-\frac 34 (4\wp(u)\wp(v)\wp(w)-g_3).\nonumber
\end{align}
\end{proposition}
\begin{proof}
  We give two proofs of these results, the first based on
  straightforward manipulations of (\ref{A1}) and (\ref{FS}), and
  the second based on a pole argument.  As we use both techniques in
  the genus 2 case, we give some detail here for completeness, and to
  aid understanding.

{\bf First Proof}.
In (\ref{A1}), put $v=\zeta u$ and use (\ref{CM1a}) to get
\[
\frac{\sigma((1-\zeta)u)}{\sigma(u)^3}=(1-\zeta)\wp(u).
\]
Next put $w=\zeta u$ in (\ref{FS}) and use the above result, the fact
that $\sigma$ is an odd function of its argument, and
(\ref{CM1a}) to give (\ref{A1e3}).
 
Now consider (\ref{A1e3}).  
Firstly, we make use of (\ref{A1}) by
taking $(u,v)$ as $(v,w)$, $(\zeta v, \zeta^2 w)$, $(\zeta^2 v,\zeta
w)$ in turn.  Multiplying all three versions together, we get
\begin{equation}\label{FS3}
  \begin{aligned}
  &\frac{\prod_{j=0}^{j=2}\sigma(u-\zeta^{2j}w)\sigma(\zeta^{j}v-\zeta^{2j}w)
    \sigma(u-\zeta^{j}v)\sigma(u+\zeta^{j}v+\zeta^{2j}w)}
  {\prod_{j=0}^{j=2}\sigma(u)\sigma(\zeta^{j}v)^3\sigma(\zeta^{2j}w)^3}\\
  =&-\frac18\left|\,
  \begin{matrix}
     1 & \wp(u) & \wp'(u)\\
     1 & \wp(v) & \wp'(v)\\
     1 & \wp(w) & \wp'(w)
  \end{matrix}\,
  \right|\left|\,
  \begin{matrix}
     1 & \wp(u) & \wp'(u)\\
     1 & \wp(\zeta v) & \wp'(\zeta v)\\
     1 & \wp(\zeta^2 w) & \wp'(\zeta^2 w)
  \end{matrix}\,
  \right|\left|\,
  \begin{matrix}
     1 & \wp(u) & \wp'(u)\\
     1 & \wp(\zeta^2 v) & \wp'(\zeta^2 v)\\
     1 & \wp(\zeta w) & \wp'(\zeta w)
  \end{matrix}\,
  \right| 
  \end{aligned}
  \end{equation}
Now note the denominator of the l.h.s.\ simplifies using (\ref{CM1a}) to
\[
\sigma(u)^9\sigma(v)^9\sigma(w)^9.
\]
Consider now the r.h.s.  Multiply this out, simplify using
(\ref{CM1a}), then replace all occurrences of $\wp(\cdot)^3$ with
$\frac14(\wp'(\cdot)^2+g_3)$.  Then factor (Maple is useful for this
calculation!).  The result is
\begin{eqnarray}
&&-\frac{1}{32}(\wp'(u)-\wp'(w)) (\wp'(v)-\wp'(w)) (\wp'(v)-\wp'(u))\times
\nonumber\\
&& \qquad \times (\wp'(v)
  \wp'(u)+\wp'(v) \wp'(w)+\wp'(w) \wp'(u)+3 g_3-12 \wp(u) \wp(v)
  \wp(w))\nonumber
\end{eqnarray}
Now apply (\ref{A1e2}) with the minus sign and with $(u,v)=
(u,v),(u,w),(v,w)$ in turn to the numerator of the l.h.s.\ of
(\ref{FS3}).  Finally, cancelling common factors, we have
(\ref{A1e3}).

{\bf Second Proof}.
Both sides of (\ref{A1e2}) are elliptic functions, by (\ref{L2.14.1}).  
Fixing $u$ and regarding both sides of (\ref{A1e2}) as a function of $v$, 
we see both sides have the same poles and zeroes with the same order at
\[
  v=0 \ \mbox{of order $-3$}, 
  \quad u \ \mbox{of order $1$}, 
  \quad \zeta u \ \mbox{of order $1$}, 
  \quad \zeta^2u \ \mbox{of order $1$},
  \]
and with no poles or zeroes elsewhere, because of (\ref{CM1a}). 
Hence the two sides coincide up to a non-zero multiplicative constant. 
Looking at the coefficients of Laurent expansion with respect to $v$, 
we see the two sides are equal. 

Now consider (\ref{A1e3}).  Recall that, in this case, $J$ is
isomorphic to $C$, and the space of functions on $J$ having a pole
only at $0$ of order at most $n$ is given by
  \[
  \Gamma(J,\mathcal{O}(n\cdot \circ )) =\begin{cases}
    \mathbb{C}                      & \mbox{if $n=0$ or $1$},\\
    \mathbb{C}\oplus\mathbb{C}\,\wp(u)   & \mbox{if $n=2$}, \\
    \Gamma(J,\mathcal{O}((n-2)\cdot \circ))\oplus\mathbb{C}\,\wp^{(n-1)}(u) &
    \mbox{if $n\geq 3$},
  \end{cases}
  \]
where we denote $\Theta^{[0]}$ by $\circ$, that is the origin of $J$.
The function  $\sigma(u)$  is expanded as
  \[
  \sigma(u)=u-\tfrac1{120}\,({\mu_3}^2+\mu_6)\,u^7+O(u^{13}).
  \]
  The left hand side of (\ref{A1e3}) is invariant under
  $u\leftrightarrow \zeta u$, $v\leftrightarrow \zeta v$,
  $w\leftrightarrow \zeta w$, and all exchanges of $u$, $v$, and $w$.
  It is an even function under $u\leftrightarrow -u$, $v\leftrightarrow
  -v$, $w\leftrightarrow -w$ simultaneously.  Moreover, it is of
  homogeneous weight $-6$.  Hence, it must be of the form
  \[
  a\,\big(\wp'(u)\wp'(v)+\wp'(u)\wp'(w)+\wp'(v)\wp'(w)\big)
  +b\,\wp(u)\wp(v)\wp(w)+c\,\mu_6,
  \]
  where $a$, $b$, $c$ are constants independent of $g_3$.  Then, by
  using first few terms of the power series expansion with respect
  to $u$ or $v$, and by balancing the two sides, we determine these
  coefficients to obtain {\rm (\ref{A1e2})}.
\end{proof}

\begin{remark}\label{remarks}
In the ``rational'' case, $\mu_3=\mu_6=0$, $\sigma(u)=u$,
$\wp(u)=1/u^2$, the formula {\rm (\ref{A1e3})} becomes the
well-known identity
  \[
  (a + b + c)(a + \zeta b + \zeta^2 c)(a + \zeta^2 b + \zeta c)
   = a^3 + b^3 + c^3 - 3abc.
  \]
\end{remark}
\begin{remark}
  Formula (\ref{A1e3}) turns up as a special case in the study of
  exceptional completely decomposable quasi-linear (CDQL) webs
  globally defined on compact complex surfaces \cite{pp10}.
\end{remark}

\subsection{Lemniscate case {\rm (}the $(2,2[1]+1)$-curve{\rm )}}
\noindent
For the curve  $y^2=x^3+\mu_4x$, we have
  \begin{equation}\label{5.02}
  \begin{aligned}
  \wp'(u)&=2y\\ 
  \wp(u)&=x
  \end{aligned}\bigg\} \ \ 
  \mbox{if}\ \ u=\int_{\infty}^{(x,y)}\frac{dx}{2y+\mu_3}, 
  \end{equation}
and
  \begin{equation}\label{5.03}
  (\wp')^2 = 4 \wp^3 + 4\mu_4\wp. 
  \end{equation}
As usual by putting  $\mu_4=-g_2/4$, we rewrite (\ref{5.03}) as
\begin{equation}\label{Wl}
  (\wp')^2 = 4 \wp^3 - g_2\wp. 
\end{equation}
This is usually called the {\em lemniscate} case (see Abramowitz \&
Stegun \cite{as72}).  By (\ref{CM}), we see that, for the $\wp$
satisfying {\rm (\ref{Wl})},
\begin{equation}\label{CM1b}
  \sigma({\i}u) = \i\sigma(u), \quad \wp({\i}u) = - \wp(u), 
  \quad \wp'({\i}u) = \i\wp'(u). 
\end{equation}
In this case, we have from (\ref{A1}) with $v \rightarrow iv$ and
(\ref{CM1b}) that
  \begin{equation}\label{Ai1}
   -\frac{\sigma(u+iv)\sigma(u-iv)}{\sigma(u)^2\sigma(v)^2}=\wp(u)+\wp(v).
  \end{equation}
Generalizing (\ref{A1}) and this, the main results for the lemniscate
case are the following addition formulae:
\begin{proposition}
\label{Tl1}
\begin{equation}
\begin{aligned}\label{Ai2}
& \frac{ \sigma(u+v+w)\sigma(u+v-w)\sigma(u-v+w)\sigma(u-v-w) }
       { \sigma(u)^4\sigma(v)^4\sigma(w)^4 }\\
&\qquad = \frac{1}{16} g_2^2 
 + \frac12 g_2 \left(\wp(v)\wp(w)+\wp(u)\wp(w)+\wp(u)\wp(v)\right)\\
&\qquad +\wp(u)^2\wp(v)^2 +\wp(u)^2\wp(w)^2+\wp(w)^2\wp(v)^2\\
&\qquad -2\wp(u)\wp(v)\wp(w)(\wp(u)+\wp(v)+\wp(w))
\equiv E_0(u,v,w).
\end{aligned}
\end{equation}
\begin{equation}
\begin{aligned}
  & \frac{ \sigma(u+iv+w)\sigma(u+iv-w)\sigma(u-iv+w)\sigma(u-iv-w) }
  { \sigma(u)^4\sigma(v)^4\sigma(w)^4 } \\
  &\qquad = \frac{1}{16} g_2^2
    + \frac12 g_2 (\wp(u)\wp(w)- \wp(v)\wp(w) -\wp(u)\wp(v))\\
  &\qquad +\wp(u)^2\wp(v)^2 +\wp(u)^2\wp(w)^2+\wp(w)^2\wp(v)^2\\
  &\qquad +2\wp(v)\wp(u)\wp(w)(\wp(w)-\wp(v)+\wp(u))\equiv E_1(u,w;v).
\end{aligned}
\end{equation}
By symmetry we have two further formulae under the transformations
$u\rightarrow iu$ and $w\rightarrow iw$, and finally we
have the $16$-term formula
\[
   \frac{ \prod_{n,m=0,1,2,3}\sigma(u+i^n v+i^m w) } {
    \sigma(u)^{16}\sigma(v)^{16}\sigma(w)^{16} } =
  E_0(u,v,w)E_1(u,v;w)E_1(u,w;v)E_1(v,w;u).
\]
\end{proposition}
\begin{proof}
The 4-term formulae are constructed from
products of the relation (\ref{FS}) in the
same way as in \S \ref{s5e}.  Similar relations for other
permutations of terms can also be constructed.
\end{proof}

\begin{remark}
The addition formulae in this section can be proved by another method, 
as in {\rm \S \ref{equipentamic}} below.
\end{remark}
\section{Genus Two}\label{s6}
\noindent
In this section we treat curves of genus two. 
So, \(am=4\) or \(5\).  
\subsection{Basis of spaces of Abelian functions}
\noindent
Using the functions in (\ref{wp_fcts}),
we denote  
  \begin{equation}\label{Delta}
  \begin{aligned}
  \Delta&=\det[\wp_{ij}]=\wp_{11}\wp_{22}-{\wp_{12}}^2, \\
  \Delta_j&=\frac{\partial}{\partial u_j}\Delta, \quad
  \Delta_{ij}=\frac{\partial^2}{\partial u_j\partial u_i}\Delta, \quad
  \mbox{etc.}
  \end{aligned}
  \end{equation}
\begin{lemma}
\label{basis2}
Let $n\geq 2$ be an integer.  The space
$\varGamma(J,\mathcal{O}(n\Theta^{[1]})$ of the functions having no
pole outside $\Theta^{[1]}$ and at most of order $n$ on $\Theta^{[1]}$
is given recursively by
\begin{equation*}
\begin{aligned}
\varGamma(J, \mathcal{O}(2\Theta^{[1]}))&
  =     \mathbb{C}\,1
   \oplus\mathbb{C} \wp_{11}
   \oplus\mathbb{C} \wp_{12}
   \oplus\mathbb{C} \wp_{22},\\
    \varGamma(J,\mathcal{O}\big((n+1)\Theta^{[1]}\big))&= 
     \tfrac{\partial}{\partial u_1}\varGamma(J,\mathcal{O}(n\Theta^{[1]}))
 \cup\tfrac{\partial}{\partial u_2}\varGamma(J,\mathcal{O}(n\Theta^{[1]})). 
  \end{aligned}
 \end{equation*}
In particular, 
\begin{equation*}
\begin{aligned}
\varGamma(J,\mathcal{O}(3\Theta^{[1]}))&= 
      \varGamma(J,\mathcal{O}(2\Theta^{[2]}))
    \oplus\mathbb{C} \wp_{111}
   \oplus\mathbb{C} \wp_{112}
   \oplus\mathbb{C} \wp_{122}
   \oplus\mathbb{C} \wp_{222}
   \oplus\mathbb{C}\,\Delta,\\
\varGamma(J,\mathcal{O}(4\Theta^{[1]}))&= 
      \varGamma(J,\mathcal{O}(3\Theta^{[1]}))
   \oplus\mathbb{C} \wp_{1111}
   \oplus\mathbb{C} \wp_{1112}
   \oplus\mathbb{C} \wp_{1122}
   \oplus\mathbb{C} \wp_{1222}
   \oplus\mathbb{C} \wp_{2222}\\
  &\quad
   \oplus\mathbb{C} \Delta_1
   \oplus\mathbb{C} \Delta_2,\\
\varGamma(J,\mathcal{O}(5\Theta^{[1]}))&= 
      \varGamma(J,\mathcal{O}(4\Theta^{[1]}))
   \oplus\mathbb{C} \wp_{11111}
   \oplus\mathbb{C} \wp_{11112}
   \oplus\mathbb{C} \wp_{11122}
   \oplus\mathbb{C} \wp_{11222}\\
  &\quad
   \oplus\mathbb{C} \wp_{12222}
   \oplus\mathbb{C} \wp_{22222}
   \oplus\mathbb{C} \Delta_{11}
   \oplus\mathbb{C} \Delta_{12}
   \oplus\mathbb{C} \Delta_{22}.\\
  \end{aligned}
 \end{equation*}
\end{lemma}
For the convenience to the reader, we list their weight below:
{\small
\[
\begin{tabular}{c|ccccccccccc}
\hline
function & \(\wp_{11}\) & \(\wp_{12}\) & \(\wp_{22}\) & \(\wp_{111}\) & \(\wp_{112}\) & \(\wp_{122}\) & \(\wp_{222}\) & \(\Delta\) & \(\Delta_1\) & \(\Delta_2\) & \(\wp_{1111}\) \\
\hline
weight   &    \(-6\)    &     \(-4\)   &     \(-2\)   &    \(-9\)     &   \(-7\)      &    \(-5\)     &     \(-3\)    &   \(-8\)   &  \(-11\)     &  \(-9\)      &    \(-12\)     \\
\hline
\end{tabular}
\]
\[
\begin{tabular}{cccccccccc}
\hline
\(\wp_{1112}\) & \(\wp_{1122}\) & \(\wp_{1222}\) & \(\wp_{2222}\) & \(\wp_{11111}\) & \(\wp_{11112}\) & \(\wp_{11122}\) & \(\wp_{11222}\) & \(\wp_{12222}\) & \(\wp_{22222}\) \\
\hline
    \(-10\)     &   \(-8\)       &     \(-6\)     &   \(-4\)       &    \(-15\)      &    \(-13\)      &   \(-11\)       &     \(-9\)      &   \(-7\)        &      \(-5\)      \\
\hline
\end{tabular}
\]
}
\begin{proof} 
  This is shown in Cho \& Nakayashiki \cite{cn06}, see especially the
  example for $g=2$ in Section 9 of that paper.
\end{proof}

\subsection{Equipentamic case} 
\noindent
We propose the name ``Equipentamic'' for the \((2,5[1])\)-curve
  \[
  f(x,y)=y^2-(x^5+\mu_{10}).
  \]
In this case we have
  \[
  \begin{aligned}
  \wp_{11}([-\zeta]u)&=\zeta^3\wp_{11}(u),\ 
  \wp_{12}([-\zeta]u)=\zeta^2\wp_{12}(u),\
  \wp_{22}([-\zeta]u)=\zeta\wp_{22}(u), \\
  \wp_{111}([-\zeta]u)&=-\zeta^2\wp_{111}(u),\ \cdots
  \end{aligned}
  \]
for $\zeta=\zeta_5=\exp(2\pi\i/5)$, because of \ref{CM}.  
\begin{proposition}\label{equipentamic}
We have 
   \begin{equation}
   \begin{aligned}
   &\frac{\sigma(u+v)\sigma(u+[\zeta]v)\sigma(u+[\zeta^2]v)\sigma(u+[\zeta^3]v)\sigma(u+[\zeta^4]v)}
         {\sigma(u)^5\sigma(v)^5}\\
   &=\tfrac5{18}[\wp_{122}(u)\wp_{1112}(v){+}\wp(v)\wp_{1112}(u)]
    -\tfrac5{144}[\wp_{122}(u)\Delta_{22}(v){+}\wp_{122}(v)\Delta_{22}(u)]\\
   &\ -\tfrac1{144}\,[\wp_{1112}(u)\wp_{22222}(v){+}\wp_{1112}(v)\wp_{22222}(u)]
    {-}\tfrac1{24}\,[\wp_{11111}(u){+}\wp_{11111}(v)]\\
   &\ -\tfrac1{576}\,[\Delta_{22}(u)\wp_{22222}(v){+}\Delta_{22}(v)\wp_{22222}(u)]
    +\tfrac{1}{24}\,\mu_{10}\,[\wp_{22222}(u)+\wp_{22222}(v)],
   \end{aligned}\label{eqi_p1}
   \end{equation}
where $\zeta=\exp(2\pi\i/5)$, 
and $[\zeta]v=[\zeta](v_1,v_2)=(\zeta v_1,\zeta^2 v_2)$.
Alternatively, the r.h.s.\ of the relation above can be written as 
 \begin{equation}
   \begin{aligned}
&= \tfrac14\,\wp_{22}(u) \wp_{222}(u)\left(
    \wp_{22}(v)\wp_{12}(v)^{2} -\wp_{11}(v) \wp_{22}(v)^{2}-4\,
    \wp_{12}(v)\wp_{11}(v) \right)\\
  &\quad+\tfrac12\,\wp_{122}(u) \left( \wp_{12}(v)\wp_{11}(v)+\wp_{22}(v) 
      \wp_{12}(v)^{2}-\wp_{11}(v) \wp_{22}(v)^{2} \right)
  \\
  &\quad -\tfrac12\,\wp_{11}(u)\wp_{111}(u)
  +\mu_{10}\wp_{22}(u)\wp_{222}(u) + (u \Leftrightarrow v).
   \end{aligned}\label{eqi_p2}
   \end{equation}
\end{proposition}
\begin{proof}
The left hand side is an Abelian function of weight $-15$ as a
function of $u$ (resp. $v$) because of (\ref{L2.14.1}), 
and has poles only along $\Theta^{[1]}$ of order $5$
by (\ref{L2.14.3}), (\ref{L2.14.4}). 
It is also invariant under $u\leftrightarrow [\zeta]u$,
$v\leftrightarrow [\zeta]v$, and $u\leftrightarrow v$.  
According to Lemma \ref{L2.15}, 
it must be a linear combination of homogeneous weight \(-15\) 
of terms of the form
  \begin{equation}\label{term}
  c_j\,\big(X_j(u)Y_j(v)+X_j(v)Y_j(u)\big),  
  \end{equation}
where \(X_j(u)\) and \(Y_j(u)\) are members of the list just below
the Lemma \ref{basis2}, with coefficients \(c_j\) being polynomial
of \(\mu_{10}\) over the rationals.  So, there are \(6\) possible
terms of (\ref{term}), namely those in appearing in the right hand
side of (\ref{eqi_p1}).  Its coefficients follows from expanding
both sides in power series in $u$ and $v$ with first several terms
after multiplying \(\sigma(u)^5\sigma(v)^5\) to both sides, by
computer calculation using Maple. To prove (\ref{eqi_p1}), we use
the known expansions of the 4-index $\wp_{ijkl}$ relations in this case
  \begin{equation*}
  \begin{aligned}
  \wp_{2222} &= 6 \wp_{22}^2 + 4 \wp_{12},\qquad
  \wp_{1222} =  6 \wp_{22} \wp_{12}-2 \wp_{11}, \\
  \wp_{1122} &= 4 \wp_{12}^2 + 2 \wp_{11} \wp_{22},\qquad
  \wp_{1112} = 6 \wp_{12} \wp_{11} - 4 \mu_{10}\\
  \wp_{1111} &= 6 \wp_{11}^2-12 \wp_{22} \mu_{10} ,
  \end{aligned}
  \end{equation*}
together with the derivatives of these equations with respect to the  $u_i$.  
In addition we use the known quadratic 3-index relations
 \begin{equation*}
  \begin{aligned}
  \wp_{222}^2 &= 4 \wp_{22}^3 +4 \wp_{11} +4 \wp_{22} \wp_{12},  \\
  \wp_{122} \wp_{222} &= -2 \wp_{22} \wp_{11} +4 \wp_{22}^2 \wp_{12} +2 \wp_{12}^2 ,\\
  \dots &= \dots
  \end{aligned}
  \end{equation*}
These substitutions lead eventually to (\ref{eqi_p2}).
\end{proof}

\begin{remark}
  Other types of addition formulae exist, for example, for \(u\),
  \(v\), \(w\in \mathbb{C}^2\), {\small
  \begin{equation}
  \frac{\sigma(u{+}v{+}w)
        \sigma(u{+}[\zeta]v{+}[\zeta]^2w)
        \sigma(u{+}[\zeta]^2v{+}[\zeta]^4w)
        \sigma(u{+}[\zeta]^3v{+}[\zeta]w)
        \sigma(u{+}[\zeta]^4v{+}[\zeta]^3w)}
       {\sigma(u)^5\sigma(v)^5\sigma(w)^5}
  \end{equation}
}is expressed in terms of \(\wp\)-functions and their derivatives. 
But it would need a big calculation to get the explicit expression. 
\end{remark}

\subsection{The $(2,2[2]{+}1)$-curve}\label{s7.3}
\noindent
Here we treat the $(2,2[2]{+}1)$-curve \(C\) given by
\begin{equation}\label{7.05}
    f(x,y)=y^2-(x^5+\mu_{4}x^3 +\mu_{8}x).
\end{equation}
The result here is not so interesting because it is essentially a
product of two of (\ref{A2}) in the Introduction.  For completeness we
describe the result here in compressed form.  This curve \(C\) has the
automorphism
\[ 
  {[\i]} : (x, y) \mapsto (-x, {\i}y), \quad 
  \Big(\frac{dx}{2y},\frac{xdx}{2y}\Big) \mapsto 
  \Big(\i\frac{dx}{2y},-{\i}\frac{xdx}{2y}\Big)
  \]
We see that
\begin{equation}\label{7.1}
  \begin{aligned}
    {[\i]}^2&=[-1]: (x, y) \mapsto (x,-y),\\
    [\i]&:(u_1,u_2) \mapsto ({\i}u_1,-{\i}u_2), \quad
    [\i]^2(u_1,u_2)\mapsto (-u_1,-u_2),\\
    \wp_{11}&([\i]u)=-\wp(u), \ \ \wp_{12}([\i]u)=\wp_{12}(u), \\
    \wp_{22}&([\i]u)=-\wp(u), \ \ \wp_{111}([\i]u)=\i\wp_{12}(u), \ \
    \mbox{etc..}
  \end{aligned}
\end{equation}
We trivially have the following formula.
\begin{equation}\label{7.9}
  \begin{aligned}
  &\frac{\sigma(u+v)\sigma(u+[\i]v)\sigma(u+[\i]^2v)\sigma(u+[\i]^3v)}
  {\sigma(u)^4\sigma(v)^4} \\
    &\qquad\qquad =\left(\wp_{11}(u)-\wp_{11}(v)+\wp_{12}(u)\wp_{22}(v)
      -\wp_{22}(u)\wp_{12}(v) \right)\\
    &\qquad\qquad\qquad \times \left(\wp_{11}(u)+\wp_{11}(v)-\wp_{12}(u)\wp_{22}(v)
      -\wp_{12}(v)\wp_{22}(u) \right).
  \end{aligned}
\end{equation}
This is not new because \([\i]^2\) is no other than the standard
involution \(u\to-u\).  We remark further on this in the final
section.

\section{Higher genus and non-hyperelliptic curves}
\label{s7}
\noindent
These new results described above for genus one and two were inspired
by results for the trigonal genus three case \cite{eemop07}.
In that paper we derived a three-term two-variable addition formula
which generalises (\ref{A1e2}) to the purely trigonal case of genus
three
\begin{equation}\label{c34}
f(x,y)=y^3-(x^4+\mu_3x^3+\mu_6x^2+\mu_9x+\mu_{12}).
\end{equation}
We have also recently proved the existence of a similar three-term
two-variable addition formula for the purely trigonal case of genus
four \cite{bego08}
\begin{equation}\label{c35}
f(x,y)=y^3-(x^5+\mu_3 x^4+\mu_6x^3 +\dots +\mu_{15}).
\end{equation}
In addition we showed in this paper the existence of three-term
three-variable addition formula for (\ref{c34}) and (\ref{c35}) which
generalise (\ref{A1e3}).  

\subsection{The three-term three-variable addition formula for the (3,4)
  purely trigonal curve}
The formula in this case generalises the three-variable formula given
in (\ref{A1e3}) for the genus 1 case.  This type of formula is
expected to be quite complicated, as the family of members in the
corresponding natural basis in each case is large.

For example in the case of (\ref{c34}), we have an addition formula of the type
\begin{equation}
\begin{aligned}
  &\frac{\sigma(u+v+w)\sigma(u+[\zeta]v+[\zeta^2]w)
    \sigma(u+[\zeta^2]v+[\zeta]w)}
  {\sigma(u)^3\sigma(v)^3\sigma(w)^3}\\
  &\qquad\qquad\qquad\qquad = \sum_{i=1}^{27} \sum_{j=1}^{27} \sum_{k=1}^{27}c_{ijk}
  U_i(u)V_j(v)W_k(w).
\end{aligned}\label{3term3}
\end{equation}
where the functions $U_i$, $V_j$ $W_k$ are all basis functions for the
space $\varGamma(J, \mathcal{O}(3 \Theta^{[2]}))$.  These are
enumerated in \cite{eemop07}.  
We can write the r.h.s.\ as
\[
C_{30}+C_{27}+ \dots C_{0},
\]
where $C_n$ has weight $-n$ in $u,v,w$ combined, and weight $n-30$ in
the $\lambda_i, i=3,\dots,0$.  Both the l.h.s.\ and the r.h.s.\ are
symmetric under all permutations in $(u,v,w)$.  So far only the terms
from $C_{30}$ up to $C_{18}$ have been calculated, and a full
description of this formula will be given elsewhere.  To illustrate
the complexity we give here the formula for $C_{30}$

\begin{align*}
  C_{30} &= \tfrac{1}{6} \wp_{13}(u) \partial_3 Q_{1333}(v)
  \wp_{111}(w) +\tfrac{1}{48} \partial_3 Q_{1333}(u) \partial_3
  Q_{1333}(v) \wp^{[22]}(w) \\
  & \qquad -\tfrac{3}{2} \wp_{12}(u) \wp_{11}(v) \wp^{[23]}(w)
  -\wp_{13}(u) \wp^{[22]}(v) \wp^{[22]}(w)\\
  & \qquad +\tfrac{1}{8} \wp^{[12]}(u) \partial_3 Q_{1333}(v)
  \wp_{112}(w) -\tfrac{1}{2} \wp_{11}(u) \wp_{11}(v) \wp_{11}(w) \\
  & \qquad -\tfrac{3}{16} \wp_{222}(u) \wp^{[22]}(v) \wp_{112}(w)
  -\tfrac{3}{8} \wp_{22}(u) \wp^{[23]}(v) \wp^{[23]}(w) \\
  & \qquad -\tfrac{1}{8} \wp^{[11]}(u) \wp^{[22]}(v) \wp^{[22]}(w)
  +\tfrac{3}{16} \wp_{122}(u) \wp_{122}(v) \wp^{[22]}(w) \\
  & \qquad -\tfrac{1}{8} \wp^{[13]}(u) \wp^{[13]}(v) \wp^{[13]}(w)
  +\tfrac{3}{8} \wp_{123}(u) \wp_{123}(v) \wp^{[33]}(w) \\
  & \qquad -\tfrac{1}{8} Q_{1333}(u) Q_{1333}(v) \wp^{[33]}(w)
  -\tfrac{3}{8} \wp_{33}(u) \wp^{[33]}(v) \wp^{[33]}(w) \\
  & \qquad -\tfrac{1}{8} \wp_{133}(u) \partial_2 Q_{1333}(v)
  \wp^{[23]}(w) +\tfrac{1}{4} \wp_{133}(u) \wp_{11}(v) \partial_1
  Q_{1333}(w) \\
  & \qquad -\tfrac{1}{8} \wp_{133}(u) \wp^{[13]}(v) \partial_1
  Q_{1333}(w) -\tfrac{3}{8} \wp_{223}(u) \wp_{113}(v) \wp^{[33]}(w) \\
  & \qquad +\tfrac{1}{4} \wp^{[12]}(u) \wp^{[12]}(v) \wp^{[22]}(w)
  -\tfrac{1}{4} P (1,3) \wp_{122}(v) \wp_{111}(w) \\
  & \qquad +\tfrac{1}{8} \wp_{111}(u) \wp_{111}(v) +\tfrac{3}{8}
  Q_{1333}(u) \wp_{113}(v) \wp_{113}(w) \\
  & \qquad + \text{all permutations of} \quad (u,v,w).
\end{align*} 

\subsection{Other formulae for higher genus curves}
There are a number of other addition formulae waiting to be
computed in explicit form for special cases of $g>2$ hyperelliptic
curves and for other special cases of trigonal and curves with higher
gonal numbers.

For instance, let
\begin{equation*}
   f(x,y)=\left\{
  \begin{aligned}
    & y^3-(x^{am}+\mu_{3a}x^{a(m-1)}+\mu_{6a}x^{a(m-2)}+\cdots \\
    & \qquad\qquad \qquad+\mu_{3a(m-1)}x^a+\mu_{3am}),
    \quad\quad \mbox{if $\gcd(am,3)=1$},\\
    & y^3-(x^{am+1}+\mu_{3a}x^{a(m-1)+1}+\mu_{6a}x^{a(m-2)+1}+\cdots\\
    & \qquad\qquad\qquad +\mu_{3a(m-1)}x^{a+1}+\mu_{3am}x), \quad
    \mbox{if $\gcd(am,3)=3$};
  \end{aligned}\right.
\end{equation*}
and $C$ be the curve defined by $f(x,y)=0$.  This should be called
$(3,a[m])$-curve or $(3,a[m]+1)$-curve, respectively.  
Let $\zeta=\exp(2\pi\i/(3a))$.  Then $C$ has automorphisms
  \begin{eqnarray}\label{a.09}
  [\zeta]: \left\{\begin{array}{ll}
  (x,y)\mapsto (\zeta^3 x,  \zeta^a y), & \mbox{if \ $\gcd(am,3)=1$}, \\
  (x,y)\mapsto (\zeta^3 x, -\zeta^a y), & \mbox{if \ $\gcd(am,3)=3$}.
  \end{array}\right.
  \end{eqnarray}
Namely, it is acted on by $W_{3a}$, the group of $3a$-th roots of $1$.
Associated with such a curve $C$, we will have various multi-term
addition formulae.

We shall give a remark by an explicit example. 
For \((3,4[1])\)-curve  \(y^3=x^4+\mu_{12}\) and the automorphism 
  \begin{equation}\label{a.10}
  [\i]: (x,y)\mapsto (-{\i}x, y),
  \end{equation}
the action \([\i]^2\) is different from standard involution \(u\to-u\) 
on the variable space of associated Abelian functions.  
The formula that expresses
  \begin{equation}
  \frac{\sigma(u+v)\sigma(u+[\i]v)\sigma(u+[\i]^2v)\sigma(u+[\i]^3v)}
       {\sigma(u)^4\sigma(v)^4}
  \end{equation}
in terms of \(\wp\)-functions seems to be interesting. 

\section*{Acknowledgements}
\noindent
The authors are grateful for a number of useful discussions with
Victor Enolski, Emma Previato, John Gibbons and Matthew England.  Part
of the work was done whilst YO was visiting Heriot-Watt, and we would
like to acknowledge the financial support of the Edinburgh
Mathematical Society and JSPS grant-in-aid for scientific researches
No.\ 19540002 for making this visit possible.

\label{lastpage}
\end{document}